\documentclass[11pt]{article}
\usepackage{mathrsfs}
\usepackage[centertags]{amsmath}
\usepackage{amsfonts}
\usepackage{amssymb}
\usepackage{amsthm}
\usepackage{cases}
\usepackage{indentfirst}
\usepackage{epsfig}
\usepackage {graphicx}
\usepackage{color}

\date{}

\definecolor{c20}{rgb}{0.,0.7,0.}
\definecolor{c30}{rgb}{0.,0.,1.}
\definecolor{c40}{rgb}{1,0.1,0.7}
\definecolor{c50}{rgb}{1,0,0}

\def\pzx#1{\textcolor{c30}{#1}}
\def\pzx#1{#1} 

\def\ltt#1{\textcolor{c50}{#1}}
\def\ltt#1{#1} 

\newtheorem{theorem}{Theorem}[section]

\newtheorem{lemma}{Lemma}[section]

\newtheorem{remark}{Remark}[section]

\numberwithin{equation}{section}

\def\E{\operatorname*{E}}

\def\P{\operatorname*{\mathbb{P}}}
\def\C{\operatorname*{\mathbb{C}}}

\def\N{\operatorname*{\mathbb{N}}}
\allowdisplaybreaks[4]

\begin{document}

\title{ \pzx{Moments convergence} of powered \pzx{normal} extremes}
\author{Tingting Li\thanks{Corresponding author. Email address: tinalee@swu.edu.cn}
\qquad Zuoxiang Peng \\
{\small School of Mathematics and Statistics, Southwest University, Chongqing, 400715, China}}

\maketitle
\begin{quote}
{\bf Abstract.}~~In this paper, convergence for moments of powered \pzx{normal} extremes  is considered under an optimal choice of \pzx{normalizing} constants. It is shown that the rates of convergence for normalized powered normal extremes depend on the power index. However, the dependence disappears for higher-order \pzx{expansions} of moments.

\end{quote}

\section{Introduction}
 Let $M_n=\max(X_1, X_2, \dots, X_n)$ denote the partial maximum of \pzx{an} independent and identically distributed (iid) random samples from standard normal population, and let $\left|M_n\right|^t$ be the powered extremes for given power index $t>0$. Hall (1980) showed that the normalized \ltt{powered} extremes $(\left|M_n\right|^t-d_n)/c_n$ with suitable norming constants $c_n>0$ and $d_n\in \mathbb{R}$ \ltt{have} the limiting Gumbel extreme value distribution $\Lambda(x)= \exp(-e^{-x})$. That is to say, as $n\rightarrow \infty$
 \begin{eqnarray}
 \label{eq1.1}
      \mathbb{P}(\left|M_n\right|^t\leq c_n x+d_n) \rightarrow \Lambda(x), \ \ \ \ x\in \mathbb{R}.
 \end{eqnarray}
\pzx{Hall (1980) also showed that with optimal normalizing constants, the best convergence rate  for \eqref{eq1.1} is  $1/\log^2 n$ \ltt{as} $t=2$, contrary to the case of $0<t\neq 2$,  which is $1/\log n$.}

Recently, Zhou and Ling (2016) studied the higher-order expansions and established the higher-order convergence rates for \ltt{both} the cumulative distribution function (cdf) and probability density function (pdf) of normalized powered extremes. It is well known that \pzx{convergence of the cdf and the pdf may not imply the convergence of moments, see, e.g., Resnick (1987)}. Therefore, \pzx{natural questions are how about the convergence and expansions of moments \ltt{for the} normalized powered extremes, respectively}.

 Moments convergence of extremes for a sequence of iid random variables from any given cdf $F$, have been of considerable interests.
\pzx{Under some suitable conditions,} Pickands (1968) proved that moments of normalized extremes converge to the corresponding moments of the extreme value
distribution \pzx{if} $F$ is in the domain of attraction of an extreme
value distribution. Nair (1981) derived asymptotic expansions for the moments of standard normal extremes. Further, Liao et al. (2013) and Jia et al. (2015) extended \pzx{the results to that of} skew-normal distribution and general error distribution, respectively. For other \pzx{work related to moments convergence on extremes}, see, e.g. Hill and Spruill (1994), H\"{u}sler et al. (2003) and Withers and Nadarajah (2011).

\par The objective of this paper is to establish the asympotics for moments of normalized powered extremes for normal samples. \pzx{The optimal normalizing constants $c_{n}$ and $d_{n}$ given by Hall (1980) will be used throughout this paper, i.e.,
\begin{equation}\label{eq1.2}
c_n = tb_n^{t-2}-2b_n^{-2}\mathbb{I}\{t=2\}, \qquad \qquad d_n= b_n^t-2b_n^{-2}\mathbb{I}\{t=2\}
\end{equation}
with $\mathbb{I}\{\cdot\}$ denoting the indicator function, where $t>0$\ltt{, the power index,} and $b_n>0 $ \ltt{is} the solution of the equation
\begin{equation}\label{eq1.3}
 2\pi b_n^2 \exp(b_n^2)= n^2.
\end{equation}
}

\par To facilitate our description, \pzx{we cite the following results due to Zhou and Ling (2016).
With the normalizing constants $c_{n}$ and $d_{n}$ given by \eqref{eq1.2}, Zhou and Ling (2016) showed that}
\begin{eqnarray}\label{eq1.4}
\frac{d}{dx}\mathbb{P}\left(\left|M_n\right|^t\leq c_n x+ d_n \right)= \Lambda^\prime(x)\Bigg(1+b_n^{-2-2I\{t=2\}}\Big(\varpi(t,x) + b_n^{-2}\tau(t,x) + O(b_n^{-4})\Big)\Bigg)
\end{eqnarray}
holds for large $n$, where \pzx{$\Lambda^{\prime}(x)=e^{-x}\Lambda(x)$,} $\varpi(t,x)=\kappa_1(t,x) + \kappa_2(t,x)$ with
\begin{equation}\label{eq1.5}
\kappa_1(t,x)=\left\{ {{\begin{array}{*{20}c}
 { \frac{1}{2}+x+x^2,
\qquad \qquad t=2;} \hfill  \\
 { x(1-t +\frac{t-2}{2}x),
\qquad t\neq 2}  \hfill \\
\end{array} }}\right.
\end{equation}

\begin{eqnarray}\label{eq1.6}
\kappa_2(t,x) = \left\{ {{\begin{array}{*{20}c}
 { -e^{-x} \left(\frac{7}{2}+3x +x^2\right),
\qquad t=2;} \hfill  \\
 { e^{-x} \left(1+x +\frac{2-t}{2}x^2\right),
\qquad t\neq 2}  \hfill \\
\end{array} }}\right.
\end{eqnarray}
and
\begin{eqnarray}\label{eq1.7}
\tau(t,x) = \left\{ {{\begin{array}{*{20}c}
 { e^{-x}\left( \frac{43}{3} +14x +6x^2 + \frac{4}{3}x^3\right),
\qquad \qquad \qquad \qquad \qquad\qquad\qquad\qquad t=2;} \hfill  \\\\{
\begin{split}
  & xe^{-x}\left(1-t+\frac{t-2}{2}x\right)\left( 1+x +\frac{2-t}{2}x^2\right)\\
+& x^2\left( \frac{(1-t)(1-2t)}{2}+\frac{5(1-t)(t-2)}{6}x + \frac{(t-2)^2}{8}x^2\right)\\
-&e^{-x}\left(3+3x+ \frac{3}{2}x^2+ \frac{(2-t)(2t+1)}{6}x^3+\frac{(t-2)^4}{8}x^4\right)\\
+&\frac{e^{-2x}}{2}\left( 1+x + \frac{2-t}{2}x^2\right)^2,
\end{split}
\qquad \qquad \ \ t\neq 2.}  \hfill \\
\end{array} }}\right.
\end{eqnarray}

The rest of the paper is organized as follows. Section \ref{sec2} gives our main results for the asymptotic behavior of moments of normalized powered extremes. Proofs of the main results \pzx{with some auxiliary lemmas}  are deferred to Section \ref{sec3}.

\section{Main results}\label{sec2}
In this section, we provide \pzx{the} higher-order expansion of moments for normalized powered normal extremes. \pzx{For simplicity,  with normalizing constants $c_{n}$ and $d_{n}$ given by \eqref{eq1.2}, let
\begin{equation}\label{eq2.1}
m_{r,t}(n)=\E\left(\left(|M_n|^t-d_n\right)/{c_n}\right)^{r}\quad \mbox{and}\quad m_{r}=\int_{-\infty}^{\infty}x^{r}d\Lambda(x),\quad r\in {\N}^{+}
\end{equation}
denote the $r$-th moment of normalized powered normal extremes with power index $t>0$ and $m_r$ be the $r$-th moment of Gumbel extreme value distribution, respectively. }

\begin{theorem} \label{thm2.1}
\pzx{For $m_{r,t}(n)$ and $m_{r}$ given by \eqref{eq2.1}, the following results hold.}
\begin{itemize}
\item[(1).]
For $0<t\neq 2$,
\begin{eqnarray}\label{eq2.2}
&&\lim_{n\to\infty}b_n^2\left[b_n^2\left(m_{r,t}(n)-m_r\right)+r\left(m_{r-1}+m_r-\frac{t-2}{2}m_{r+1}\right)\right]\nonumber\\
&=& \frac{5}{2}r m_{r-1}
 +\left( \left(t+1\right) \left(r+1\right)-\frac{5}{2}\right)m_{r}
 + \left( \left(r+1\right)t -1 \right)m_{r+1} \nonumber\\
 &&+ \left( \frac{5\left(2-t\right) \left(t-1\right) \left(r+3\right) }{6} + \frac{2t^2-5t+1}{2}\right) m_{r+2} \nonumber\\
 &&+ \frac{(t-2)^2(r+4)}{4}m_{r+3}
 -\frac{(t-2)^2}{8}m_{r+4}.
\end{eqnarray}
\item[(2).] For $t=2$,
\begin{eqnarray}\label{eq2.3}
&&\lim_{n\to\infty}b_n^2\left[b_n^4\left(m_{r,2}(n)-m_r\right)+r\left(\frac{7}{2}m_{r-1}+3m_r+m_{r+1}\right)\right]\nonumber\\
&= & -\frac{43}{3}r m_{r-1}
 +\left( \frac{1}{3}-14r\right)m_{r}
 + \left( 2-6r \right)m_{r+1} \nonumber\\
 &&+ \left( 2-\frac{4}{3}r \right) m_{r+2} + \frac{4}{3}m_{r+3}.
\end{eqnarray}
\end{itemize}
\end{theorem}

\begin{remark}\label{remark2.1}
Noting that $b_n^2\sim 2\log n $ from \eqref{eq1.3}, we know that for large $n$
\begin{equation*}
m_{r,t}(n)- m_r \sim -\frac{r(m_{r-1}+m_r-\frac{t-2}{2}m_{r+1})}{2\log n}
\end{equation*}
as $0<t\neq 2$, and for $t=2$
\begin{equation*}
m_{r,2}(n)- m_r \sim -\frac{r\left(\frac{7}{2}m_{r-1}+3m_r+m_{r+1}\right)}{4\log^2 n}.
\end{equation*}
\end{remark}

\begin{remark}\label{remark2.2}
Similar to the findings of Hall (1980) and Zhou and Ling (2016) \ltt{respectively on the convergence rates of the cdf and the pdf of normalized extremes}, Remark \ref{remark2.1} implies that moments of normalized powered normal extremes converge to the corresponding moments of the limiting extreme value distribution, and the rates of convergence depend on the power index $t$. \pzx{Precisely}, the moments of normalized squares of normal extremes $M_n^2$ converge at a rate of $1/\log^2 n$ while for $0<t\neq 2$,  the rates of convergence become $1/\log n$. However, Theorem \ref{thm2.1} shows that this difference disappears for higher-order expansion of moments.
\end{remark}

\section{Proofs}\label{sec3}
\par In this section, we provide the proofs of our main results \ltt{in} Theorem \ltt{\ref{thm2.1}}. We shall only prove \ltt{\eqref{eq2.2}} \ltt{for the case of $0< t \neq 2$ since by using similar techniques \eqref{eq2.3} can be established as $t=2$} . To facilitate our proofs, some auxiliary lemmas are \pzx{needed}. In the sequel, the symbol $\C$ will denote a generic positive constant depending only on $r, t$, and $s(x)$ will denote a generic polynomial on $x$ independent of $n$. \pzx{For simplicity, let
\begin{equation}\label{eq3.00}
z_n= (c_n x+ d_n)^{1/t}\quad \mbox{and}\quad \alpha_n = b^2_n/4(4+t)
\end{equation}
with $x\ge -d_{n}/c_{n}$ and $t>0$, where constants $c_{n}$, $d_{n}$ and $b_{n}$ are respectively given by \eqref{eq1.2} and \eqref{eq1.3}. \ltt{Denote}
\begin{equation}\label{eq3.0}
C_{n}(x)= n\phi(z_n) \frac{dz_n}{dx}
\end{equation}
with $\phi(\cdot)$ \ltt{being} the pdf of $N(0,1)$.}

\begin{lemma}\label{lemma3.1}
If $|x|\leq \alpha_n$, for large $n$,
\begin{equation*}
 b_n^{4} \left|\exp\left(\frac{1}{2}b_n^2\left(1-\left(1+tb_n^{-2}x\right)^{\frac{2}{t}}\right)\right)- \left(1+\frac{t-2}{2}b_n^{-2}x^2\right)e^{-x}\right|\leq s(x)\exp\left(-x+ \frac{1}{4}|x|\right),
\end{equation*}
where $s(x)\geq 0$ is a polynomial on $x$ \ltt{independent of $n$}.
\end{lemma}
\noindent\textbf{Proof.} \pzx{It is obviously true for the case of $t=2$ with $s(x)=0$. The rest is for the case of $0<t\neq 2$.} By using Taylor expansion for $\left(1+tb_n^{-2}x\right)^{{2}/{t}}$ with Lagrange remainder term, for $\theta \in (0,1)$ we write
\begin{eqnarray*}
&&\exp\left(\frac{1}{2}b_n^2\left(1-\left(1+tb_n^{-2}x\right)^{\frac{2}{t}}\right)\right)\\
&&=\exp\left(\frac{1}{2}b_n^2\left(1- \left(1+ 2b_n^{-2}x + \left(2-t\right)b_n^{-4}x^2 + \frac{2}{3} \left(2-t\right)\left(1-t\right)b_n^{-6}x^3\left(1+\theta tb_n^{-2}x\right)^{\frac{2}{t}-3}\right)\right)\right)\\
&& =\exp\left(-x + \frac{t-2}{2}b_n^{-2}x^2 -\frac{1}{3}(2-t)(1-t)b_n^{-4}x^3(1+\theta tb_n^{-2}x)^{\frac{2}{t}-3}\right).
\end{eqnarray*}
Therefore, using the \pzx{inequality} $1+x\leq e^{x}, x\in \mathbb{R}$, we have
\begin{eqnarray}\label{eq3.1}
\exp\left(\frac{1}{2}b_n^2\left(1-\left(1+tb_n^{-2}x\right)^{\frac{2}{t}}\right)\right)
&\geq &\left(1+ \frac{t-2}{2}b_n^{-2}x^2 -\frac{1}{3}(2-t)(1-t)b_n^{-4}x^3(1+\theta tb_n^{-2}x)^{\frac{2}{t}-3}\right)e^{-x}\nonumber\\
&\geq & \left(1+ \frac{t-2}{2}b_n^{-2}x^2 -\frac{1}{3}(2+t)(1+t)b_n^{-4}|x|^3|(1+\theta tb_n^{-2}x)^{\frac{2}{t}-3}|\right)e^{-x}\nonumber\\
&\geq &  \left(1+ \frac{t-2}{2}b_n^{-2}x^2 -(2+t)(1+t)b_n^{-4}|x|^3\right)e^{-x}.
\end{eqnarray}
The last inequality holds since $|\theta tb_n^{-2}x|\leq \frac{t}{16+4t}\leq \frac {1}{4}$ when $|x|\leq \alpha_n $, \pzx{and
\begin{equation*}
\Big|\left(1+\theta tb_n^{-2}x \right)^{\frac{2}{t}-3}\Big|\leq\max\left(\left(1+ \frac{t}{16+4t}\right )^{\frac{2}{t}},\left(1-\frac{1}{4} \right)^{\frac{2}{t}-3}\right) <3.
\end{equation*}
}On the other hand, by using the inequality $\exp(x) \leq 1+x + \frac{1}{2}x^2\exp(|x|), x\in \mathbb{R}$, \pzx{for large $n$ we have}
\begin{eqnarray}\label{eq3.2}
&&\exp\left(\frac{1}{2}b_n^2\left(1-\left(1+tb_n^{-2}x\right)^{\frac{2}{t}}\right)\right)\nonumber\\
&\leq& e^{-x}\left[ 1+ \frac{t-2}{2}b_n^{-2}x^2 -\frac{1}{3}(2-t)(1-t)b_n^{-4}x^3(1+\theta tb_n^{-2}x)^{\frac{2}{t}-3}\right.\nonumber\\
&&+ \frac{1}{2} \left(\frac{t-2}{2}b_n^{-2}x^2 -\frac{1}{3}(2-t)(1-t)b_n^{-4}x^3(1+\theta tb_n^{-2}x)^{\frac{2}{t}-3}\right)^2\nonumber\\
&&\left. \exp\left(\left|\frac{t-2}{2}b_n^{-2}x^2 -\frac{1}{3}(2-t)(1-t)b_n^{-4}x^3(1+\theta tb_n^{-2}x)^{\frac{2}{t}-3}\right|\right)\right]\nonumber\\
&\leq&  \left(1+ \frac{t-2}{2}b_n^{-2}x^2\right)e^{-x}+ b_n^{-4}(2+t)(1+t)|x|^3 e^{-x}\nonumber\\
&& + \frac{1}{8}b_n^{-4}(t-2)^2 x^4(1+|x|)^2\exp\left(-x + \left|\frac{t-2}{2}b_n^{-2}x^2\right| + \left|(2-t)(1-t)b_n^{-4}x^3\right| \right)\nonumber\\
&\leq & \left(1+ \frac{t-2}{2}b_n^{-2}x^2\right) e^{-x}+ b_n^{-4}s(x)\exp\left(-x+ \frac{1}{4}|x|\right),
\end{eqnarray}
where the last inequality holds since\pzx{
\begin{equation*}
\max\left(\left|\frac{t-2}{2}b_n^{-2}x^2\right|, \left|(2-t)(1-t)b_n^{-4}x^3\right|\right)\leq \frac{1}{8}|x|
\end{equation*}
if $|x|\leq \alpha_n$. The desired result follows by \eqref{eq3.1} and \eqref{eq3.2}.} \qed

\begin{lemma}\label{lemma3.2}
Let $C_{n}(x)= n\phi(z_n) \frac{dz_n}{dx}$ \pzx{given by \eqref{eq3.0}}. For
\pzx{$|x| \leq \alpha_n$ and large $n$, functions} $C_{n}(x)$,
$b_{n}^{2}(C_{n}(x)-e^{-x})$ and $b_{n}^{2}\left[
b_{n}^{2}(C_{n}(x)-e^{-x})- \kappa_1(t,x)e^{-x}  \right]$  are
bounded by $s(x)\exp(-x+\frac{1}{4}|x|)$, where $\kappa_1(t,x)$ is given by \eqref{eq1.5} \ltt{and} $s(x)\geq 0$ is a polynomial on $x$ \ltt{independent of $n$}.
\end{lemma}
\noindent\textbf{Proof.}~~Note that we only need to show that $b_{n}^{2}\left[
b_{n}^{2}(C_{n}(x)-e^{-x})- \kappa_1(t,x)e^{-x}  \right]$ can be bounded by $s(|x|)\exp\left(-x+\frac{1}{4}|x|\right)$.
\pzx{By \eqref{eq1.3} and simple calculations}, we have
\begin{eqnarray*}
&&b_{n}^{2}\Big[
b_{n}^{2}\left(C_{n}(x)-e^{-x}\right)- \kappa_1(t,x)e^{-x}  \Big]\\
&=&\left(1+tb_n^{-2}x\right)^{\frac{1}{t}-1}b_n^{4}\Big[ \exp\left(\frac{1}{2}b_n^2\left(1-(1+tb_n^{-2}x\right)^{\frac{2}{t}})\right)-\left(1+\frac{t-2}{2}b_n^{-2}x^2\right)e^{-x}\Big]\\
&&+ \Big[b_n^{4} \left (\left(1+tb_n^{-2}x\right)^{\frac{1}{t}-1}-1-(1-t)b_n^{-2}x\right)
+ \frac{t-2}{2}x^2 b_n^{2} \left(\left(1+tb_n^{-2}x\right)^{\frac{1}{t}-1}-1\right)\Big]e^{-x}.
\end{eqnarray*}
Therefore, \pzx{the desired result follows} by Lemma \ref{lemma3.1} and Taylor's expansion of $(1+tb_n^{-2}x)^{{1}/{t}-1}$.\qed

\begin{lemma}\label{lemma3.3}
Let $\Psi_n(x)= 1- \Phi(z_n)$, where $\Phi(\cdot)$ denote the cdf of $N(0,1)$.
For  $|x| \leq \alpha_n$ and large $n$, we have that $b_n^2 (e^{-x} -n\Psi_n(x))$ and  $b_n^2\{b_n^2 \left[e^{-x} -n\Psi_n(x)\right]- \kappa_2(t,x)\}$ can be
bounded by $s(x)\exp(-x+\frac{1}{4}|x|)$, where $\kappa_2(t,x)$ is given by \eqref{eq1.6} \ltt{and $s(x)\geq 0$ is a polynomial on $x$ independent of $n$}.
\end{lemma}
\noindent\textbf{Proof.}~~
For $x>0$, we have
$$ 1- \Phi(x) = x^{-1}\phi(x)(1-x^{-2}f(x)) $$
where $0<f(x)<1$. Then by Lemma 1 in Hall (1980) and Lemma \ref{lemma3.1}, we have
\begin{eqnarray}\label{eq3.3}
&&b_n^2\Big\{b_n^2 \left[e^{-x} -n\Psi_n(x)\right]- \kappa_2(t,x)\Big\}\nonumber\\
&\leq & b_n^2\Big\{b_n^2 \left[e^{-x}- nz_n^{-1}\phi(z_n)(1-z_n^{-2})\right]- \kappa_2(t,x)\Big\}\nonumber\\
& = &  b_n^2\Big\{b_n^2 \left[e^{-x}- \left(1+tb_n^{-2}x \right)^{-\frac{1}{t}}\exp\left(\frac{1}{2}b_n^2\left(1- \left(1+tb_n^{-2}x\right)^{-\frac{2}{t}}\right)\right)\right.  \nonumber\\
&& \times\left. \left(1- b_n^{-2}\left(1+tb_n^{-2}x\right)^{-\frac{2}{t}}\right)\right]- \kappa_2(t,x)\Big\}\nonumber\\
&\leq& b_n^2\Bigg\{b_n^2 \Big[e^{-x}- \left (1-b_n^{-2}x\right)\left(1- b_n^{-2}\left(1-2b_n^{-2}x+ b_n^{-4}x^2(32+2t)\right)\right)\nonumber\\
&& \times\exp\left(\frac{1}{2}b_n^2(1- (1+tb_n^{-2}x)^{-\frac{2}{t}})\right)\Big]-\kappa_2(t,x)\Bigg\}\nonumber\\
&\leq &  b_n^2\Bigg\{b_n^2 \left[e^{-x}- \left(1- b_n^{-2}(1+x)- b_n^{-4}s_1(x)\right)\exp\left[\frac{1}{2}b_n^2\left(1- \left(1+tb_n^{-2}x\right)^{-\frac{2}{t}}\right)\right]\right]-\kappa_2(t,x)\Bigg\}\nonumber\\
&=& b_n^4 e^{-x}\left[1-\left(1- b_n^{-2}(1+x)- b_n^{-4}s_1(x)\right)(1+\frac{t-2}{2}b_n^{-2} x^2)- b_n^{-2}(1+x+\frac{2-t}{2}x^2)\right]\nonumber\\
&+& b_n^4(1- b_n^{-2}(1+x)- b_n^{-4}s_1(x))\left[\exp\left(\frac{1}{2}b_n^2\left(1-\left(1+tb_n^{-2}x\right)^{\frac{2}{t}}\right)\right)- (1+\frac{t-2}{2}b_n^{-2}x^2)e^{-x}\right]\nonumber\\
&\leq& s_2(x)e^{-x}+ s_3(x)\exp\left(-x+\frac{1}{4}|x|\right)\nonumber\\
&\leq& s(x)\exp\left(-x+\frac{1}{4}|x|\right),
\end{eqnarray}
where $s_i(x)$, $i=1,2,3$ and \ltt{ $s(x)\geq 0$ } denote generic polynomials on $x$ independent of $n$.

Similarly, we can prove that
\begin{equation}\label{eq3.4}
b_n^2\Bigg\{b_n^2 \Big[e^{-x} -n\Psi_n(x)\Big]- \kappa_2(t,x)\Bigg\}\geq -s(x)\exp\left(-x+\frac{1}{4}|x|\right).
\end{equation}
\pzx{The result follows by \eqref{eq3.3} and \eqref{eq3.4}.} \qed

\begin{lemma}
\label{lemma3.4} For large n and $-d\log b_{n}< x \leq \alpha_n$ with $0< d<1$,
$$x^{r}C_n(x)b_n^{2}\Big\{b_n^{2}\left[\Phi^{n-1}(z_n)-
\Lambda( x)\right] - \kappa_2(t,x)\Lambda(x)\Big\}$$
are bounded by
integrable functions independent of $n$, where $r>0$, \pzx{where $C_{n}(x)$ is given by \eqref{eq3.0}.}
\end{lemma}
 \noindent
\textbf{Proof.}~~\pzx{Note that}
\begin{eqnarray*}
&&b_n^2\Big\{b_n^2 \left[\Phi ^{n-1}(z_n)-\Lambda(x)\right]- \kappa_2(t,x)\Lambda(x)\Big\}\\
&=& b_n^2\Lambda(x)\Big\{b_n^2 \left[\exp\left(n\log\Phi(z_n)+e^{-x}\right)-1\right]-\kappa_2(t,x)\Big\}\\
& &+ b_n^4\Lambda(x)\left(\Phi^{-1}(z_n)-1\right)\exp\left(n\log\Phi(z_n)+e^{-x}\right)\\
&\equiv & I_1 + I_2.
\end{eqnarray*}
For $I_2$, when $x>-d\log b_n$ with $0<d<1$,
$$ \Phi^{-1}(z_n)-1 \leq \frac{1-\Phi(b_n(1-dtb_n^{-2}\log b_n)^{1/t})}{\Phi(b_n(1-dtb_n^{-2}\log b_n)^{1/t})}\leq \C b_n^{-1} \exp\left(-\frac{1}{2}b_n^2 (1-dtb_n^{-2}\log b_n)^{\frac{2}{t}}\right),$$
for large $n$, where $\C$ is a positive constant.
Therefore,
\begin{equation*}
|I_2|\leq \C b_n^3\Lambda(x)\exp(b_n^d) \exp\left(-\frac{1}{2}b_n^2 (1-dtb_n^{-2}\log b_n)^{\frac{2}{t}}\right)\leq \frac{1}{2}\Lambda(x)
\end{equation*}
\pzx{for $x>-d\log b_n$ and large $n$}. Next, we deal with $I_1$. Write $n\log \Phi(z_n)= - n \Psi_n(x) -R_n(x)$ where
$$ 0< R_n(x) \leq  \frac{n\Psi_n^2(x)}{2(1-\Psi_n(x)) }.$$
Then
\begin{eqnarray*}
|I_1|  &\leq& b_n^2\Lambda(x)\left|b_n^2 \left(e^{-x}- n\Psi_n(x)\right)-\kappa_2(t,x)\right| \\
 &+& \frac{1}{2}b_n^4  \Lambda(x) (e^{-x}- n\Psi_n(x))^2 \exp( |e^{-x}- n\Psi_n(x)|) \\
 &+& b_n^4 \Lambda(x)|\exp[-R_n(x)]-1|\exp(|e^{-x}-n\Psi_n(x)|)\\
 &\equiv & II_1 + II_2 + II_3.
\end{eqnarray*}
\pzx{It follows from Lemma 2 and Lemma 3 of Hall (1980) that}\ltt{
\begin{eqnarray*}
\Big|\exp(-x) - n \Psi_n(x)\Big| \leq  \left\{ {{\begin{array}{*{20}c}
 {\C b_n^{-2}(1+|x|+x^2) \exp\left(\frac{5}{4}|x|\right)\leq 1,
\qquad -d\log b_n <x\leq0,} \hfill \\
 { \C b_n^{-2}(1+|x|+x^2) \exp\left(-\frac{3}{4}x\right)\leq \C,
\qquad \qquad 0\leq x\leq \alpha_n,} \hfill \\
\end{array} }}\right.
\end{eqnarray*}}
and
$$
\Big|\exp[-R_n(x)]-1\Big|\leq \C n^{-\frac{1}{2}},\ \ x>- d\log b_n.
$$
for large $n$. Therefore, by Lemma \ref{lemma3.2} and Lemma \ref{lemma3.3}, $x^{r}C_n(x)b_n^{2}\left\{b_n^{2}\left[\Phi^{n-1}(z_n)-
\Lambda( x)\right] - k_2(t,x)\Lambda(x)\right\}$ can be bounded by integrable function when $-d\log n <x\leq \alpha_n$. Then
we complete the proof of this lemma. \qed

\begin{lemma}\label{lemma3.5} \pzx{Let $C_{n}(x)$ be given by \eqref{eq3.0}.} For any $k>0$, we have
\begin{eqnarray*}
\lim_{n\to\infty}b_n^k\int_{-b_n^2/t}^{\infty} |x|^r C_n(x)\Phi^{n-1}(-z_n)dx = 0,\quad
\lim_{n\to\infty}b_n^k  \int_{\alpha_n}^{\infty} |x|^r C_n(x) dx = 0.
\end{eqnarray*}
\end{lemma}

\noindent\textbf{Proof.}~~ \pzx{Obviously, $\Phi(-z_n)\leq 1/2$ by noting that $z_n= b_n(1+tb_n^{-2}x)^{\frac{1}{t}}\geq 0$ since  $x\geq -b_n^2/t$. Therefore, changing} the variables in the integrals yields
\begin{eqnarray*}
&& b_n^k\int_{-b_n^2/t}^{\infty} |x|^r C_n(x)\Phi^{n-1}(-z_n)dx\\
&& \leq 2^{1-n} b_n^k \exp(\frac{b_n^2}{2})  \int_{-b_n^2/t}^{\infty} |x|^r (1+ tb_n^{-2}x)^{\frac{1}{t}-1} \exp\left(-\frac{1}{2}b_n^2\left(1+ tb_n^{-2}x\right)^{\frac{2}{t}}\right)dx \\
&& \leq \C  2^{-n} b_n^{k+2r+2}\exp(\frac{b_n^2}{2}) \sum_{j=0}^r \frac{r!}{j!(r-j)!} \int_0 ^{\infty} s^{tj} \exp(-s^2)ds \\
&& \leq \C  2^{-n} n b_n^{\ltt{k+2r+1}}\to 0
\end{eqnarray*}
\pzx{as $n\to\infty$ due to \eqref{eq1.3} and} $b_{n}^{2}\sim 2\log n$. \ltt{Now switching to the second limit, we have}
\begin{eqnarray*}
&& b_n^k\int_{\alpha_n}^{\infty} |x|^r C_n(x)dx\\
&& \leq b_n^k  \int_{\alpha_n}^{\infty} |x|^r (1+ tb_n^{-2}x)^{\frac{1}{t}-1} \exp\left( \frac{1}{2}b_n^2\left(1-(1+ tb_n^{-2}x)^{\frac{2}{t}}\right)\right)dx \\
&& \leq \C b_n^{k+2r+2} \sum_{j=0}^r \frac{r!}{j!(r-j)!}\int_{(\frac{16+5t}{16+4t})^{\frac{1}{t}}}^{\infty} s^{tj} \exp\left( \frac{1}{2}b_n^2\left(1-s^2\right)\right)ds \\
&& \leq \C b_n^{k+2r+2} \exp\left( \frac{1}{2}b_n^2\left(1- \left(\frac{16+5t}{16+4t}\right)^{\frac{1}{t}}\right)\right)
 \sum_{j=0}^r \frac{r!}{j!(r-j)!} \int_0^{\infty}  s^{tj} \exp\left(\left( 1- \left(\frac{16+5t}{16+4t}\right)^{\frac{1}{t}}\right) s \right) ds \\
&& \rightarrow 0
\end{eqnarray*}
\pzx{as $n\to\infty$.}
\qed

\noindent\textbf{Proof of Theorem \ref{thm2.1}.} \pzx{Obviously,
\begin{eqnarray*}
\P\left(|M_n|^t\leq c_n x+ d_n\right)= \Phi^n(z_n) - \Phi^n(-z_n)
\end{eqnarray*}
by noting that $z_n= (c_n x+ d_n)^{1/t}\geq 0 $ given by \eqref{eq3.00}.}

\par Let $g_{n,t}(x)$ denote the density function of normalized powered extremes $(|M_n|^t-d_n)/c_n$, then
\begin{eqnarray*}
g_{n,t}(x)= \left\{ {{\begin{array}{*{20}c}
 {n\phi(z_n) \frac{dz_n}{dx} \left(\Phi^{n-1}(z_n)+\Phi^{n-1}(-z_n)\right), \qquad x> -\frac{1}{t} b_n^2;} \hfill \\
 { 0,
\qquad \qquad \qquad \qquad \qquad \qquad \qquad \qquad \ \ \ \mbox{otherwise.}} \hfill \\
\end{array} }}\right.
\end{eqnarray*}
Then,
\begin{eqnarray}\label{eq3.7}
&&b_n^2\left[b_n^2\left(m_{r,t}(n)-m_r\right)+r\left(m_{r-1}+m_r-\frac{t-2}{2}m_{r+1}\right)\right]\nonumber\\
&=& \int_{-\infty}^{\infty} x^{r}b_n^{2}\Big[b_n^{2}(g_{n,t}(x)-
\Lambda^{\prime}( x)) -\varpi(t,x)\Lambda^{\prime}(x)\Big]dx \nonumber\\
&=&\int_{-\infty}^{-b_n^2/t} x^{r}b_n^{2}\Big[b_n^{2}\left(0-
\Lambda^{\prime}( x)\right) -\varpi(t,x)\Lambda^{\prime}(x)\Big]dx \nonumber\\
&& + \int_{-b_n^2/t}^{\infty} x^{r}b_n^{2}\Big[b_n^{2}(g_{n,t}(x)-
\Lambda^{\prime}( x)) -\varpi(t,x)\Lambda^{\prime}(x)\Big]dx\nonumber\\
&=& \int_{-b_n^2/t}^{\infty} x^{r}b_n^{2}\Big[b_n^{2}(g_{n,t}(x)-
\Lambda^{\prime}( x)) -\varpi(t,x)\Lambda^{\prime}(x)\Big]dx + o(1),
\end{eqnarray}
where $\varpi(t,x)$ is given by \eqref{eq1.4}.
Rewrite\pzx{
\begin{eqnarray*}
 &&  b_n^2\Big(g_{n,t}(x)-\Lambda^{\prime}(x)\Big)- \varpi(t,x)\Lambda^{\prime}(x)\\
&=&  C_n( x)\Big[b_n^2 \left(\Phi ^{n-1}(z_n)-\Lambda(x)\right)- \kappa_2(t,x)\Lambda(x)\Big]\\
&& +\Lambda(x)\Big[ b_n^2\left(C_n( x) - \exp(-x)\right)- \kappa_1(t,x)\exp(-x)\Big]\\
&& +  \Big[C_n( x)- \exp(-x)\Big] \kappa_2(t,x)\Lambda(x)  +  b_n^2 C_n( x) \Phi^{n-1}(-z_n)\\
&\equiv& E_n+ F_n + G_n+ H_n,
\end{eqnarray*}
}
where $C_n(x)= n\phi(z_n)\frac{dz_n}{dx}$ \pzx{ is given by \eqref{eq3.0}.
It follows from Lemma \ref{lemma3.5} that
\begin{eqnarray*}
\int_{-b_n^2/t}^{+\infty} x^{r}b_n^{2}H_n dx \rightarrow 0, \ \ \ \int_{\alpha_n}^{+\infty} x^{r}b_n^{2}\left(E_n+ F_n + G_n\right) dx \rightarrow 0.
\end{eqnarray*}
Plugging into \eqref{eq3.7} yields
\begin{eqnarray}\label{eq3.8}
&&b_n^2\left[b_n^2\left(m_{r,t}(n)-m_r\right)+r\left(m_{r-1}+m_r-\frac{t-2}{2}m_{r+1}\right)\right]\nonumber\\
&=& \int_{-b_n^2/t}^{\alpha_n} x^{r}b_n^{2}( E_n+ F_n + G_n)dx + o(1).
\end{eqnarray}
\pzx{If}  the following two facts
\begin{eqnarray}\label{eq3.9}
\int_{-b_n^2/t}^{-\alpha_n} x^{r}b_n^{2}( E_n+ F_n + G_n)dx \rightarrow 0, \ \ \  \int_{-\alpha_n}^{-d\log b_n} x^{r}b_n^{2}( E_n+ F_n + G_n)dx \rightarrow 0
\end{eqnarray}
hold as $n\to\infty$, we can obtain the desired result since}
\begin{eqnarray}\label{eq3.10}
&&b_n^2\left[b_n^2\left(m_{r,t}(n)-m_r\right)+r\left(m_{r-1}+m_r-\frac{t-2}{2}m_{r+1}\right)\right]\nonumber\\
&=&\int_{-d\log b_n}^{\alpha_n }x^{r}b_n^{2}\left( E_n+ F_n + G_n \right) dx + o(1) \nonumber \\
&\rightarrow & \int_{-\infty}^{\infty} x^r \tau(t,x) \Lambda^{\prime}(x)dx
\end{eqnarray}
as $n\to\infty$, where the last assertion holds because of Lemma \ref{lemma3.2}, Lemma \ref{lemma3.4} and the dominated convergence theorem.

Now \pzx{the remainder is to prove \eqref{eq3.9}. First} we prove that, for large $n$,
\begin{eqnarray}\label{eq3.11}
\Phi^{n-1}(z_n) \leq  \left\{ {{\begin{array}{*{20}c}
 {\exp \left( -\frac{1}{2}\exp\left( \frac{3}{4} |x|\right) \right),
\qquad -\alpha_n\le x\le-d\log b_n,} \hfill \\
 { \exp(-\frac{1}{2}\exp(\C b_n^2)),
\qquad \qquad -\frac{1}{t}b_n^{2}\le x\le -\alpha_n,} \hfill \\
\end{array} }}\right.
\end{eqnarray}
where $\C$ is a generic positive constant.
In fact, for large $n$, if $x \leq -\alpha_n$,
\begin{eqnarray*}
\Phi^{n-1}(z_n) &\leq & \Phi^{n-1}\left(b_n\left(1-\frac{t}{16+4t}\right)^{\frac{1}{t}}\right) \\
&\leq& \exp\left( - \left(n-1\right) \left(1-\Phi\left(\left(\frac{16+3t}{16+4t}\right)^{\frac{1}{t}} b_n\right)\right)   \right) \\
& \leq & \exp\left( -\frac{1}{2} \exp (\C b_n^2 ) \right),
\end{eqnarray*}
where the last inequality holds because of $1- \Phi(x)\geq x^{-1}\phi(x) (1-x^{-2})$ for $x> 0$.
\pzx{For the case of $-\alpha_n \leq x \leq -d\log b_n$ with} large $n$,
\begin{eqnarray*}
\Phi^{n-1}(z_n) &\leq & \exp\Big[- (n-1) (1-\Phi(z_n))\Big]\\
&\leq& \exp\Big[ - \left(1-1/n\right) \left(1+tb_n^{-2}x\right)^{-\frac{1}{t}} \left(1- b_n^{-2}(1+tb_n^{-2}x)^{-2/t} \right)\exp\left( \frac{1}{2} b_n^2 \left( 1-(1+tb_n^{-2}x)^{2/t}\right) \right) \Big]
\\
& \leq & \exp\Big[ -\frac{1}{2} \exp\left( \frac{1}{2} b_n^2 \left( 1-(1+tb_n^{-2}x)^{2/t}\right) \right) \Big]\\
& \leq & \exp\left( -\frac{1}{2} \exp\left( -x -4 b_n^{-2} x^2 \right) \right)\\
& \leq & \exp \Big[ -\frac{1}{2}\exp\left( \frac{3}{4} |x|\right) \Big].
\end{eqnarray*}
\pzx{Thus it follows from \eqref{eq3.11} that} for any $k>0$,
\begin{eqnarray}\label{eq3.12}
&& \left |\int_{-b_n^2/t} ^{-\alpha_n} x^r C_n(x) \Phi^{n-1}(z_n) dx \right|\nonumber \\
 &\leq & \int_{-b_n^2/t} ^{-\alpha_n} |x|^r (1+ tb_n^{-2}x)^{\frac{1}{t}-1} \exp\left(\frac{1}{2} b_n^2 (1- (1+ tb_n^{-2} x)^{\frac{2}{t}})\right)\Phi^{n-1}(z_n) dx \nonumber\\
 &\leq & \exp(\frac{1}{2} b_n^2) \exp\left( -\frac{1}{2} \exp (\C b_n^2 ) \right) \int_{-b_n^2/t} ^{-\alpha_n} |x|^r (1+ tb_n^{-2}x)^{\frac{1}{t}-1} \exp\left(-\frac{1}{2} b_n^2 (1+ tb_n^{-2} x)^{\frac{2}{t}}\right)  dx \nonumber\\
 &=& \exp(\frac{1}{2} b_n^2) \exp\left( -\frac{1}{2} \exp (\C b_n^2 ) \right)O(b_n^{2r+2}) \nonumber\\
 &= &  o(b_n^k)
\end{eqnarray}
for large $n$, and by Lemma \ref{lemma3.2} we have
\begin{eqnarray}\label{eq3.13}
 && \left|\int_{-\alpha_n} ^{-d\log b_n} x^r C_n(x)\Phi^{n-1}(z_n)dx \right|\nonumber\\
 &\leq & \int_{-\alpha_n} ^{-d\log b_n} s(x)\exp (\frac{5}{4}|x|) \exp \left( -\frac{1}{2} \exp\left( \frac{3}{4} |x|\right) \right) dx\nonumber\\
 & \leq & \exp \left( - \frac{1}{4} \exp\left( \frac{3}{4} d \log b_n\right) \right)\int_{-\alpha_n} ^{-d\log b_n} s(x)\exp (\frac{5}{4}|x|) \exp \left( - \frac{1}{4} \exp\left(  \frac{3}{4} |x|\right) \right) dx\nonumber \\
 &=& o(b_n^k).
\end{eqnarray}
Similarly, we can also prove
\begin{eqnarray}\label{eq3.14}
\int_{-b_n^2/t} ^{-d\log b_n} x^r C_n(x)\exp(-x)\Lambda(x)dx= o(b_n^{k}), \ \ \int_{-b_n^2/t} ^{-d\log b_n} x^r \exp(-2x)\Lambda(x)dx= o(b_n^{k}).
\end{eqnarray}
Thus it is easy to show that \eqref{eq3.9} holds by \eqref{eq3.12}-\eqref{eq3.14}.  The proof is complete. \qed

\end{document}